\documentclass[12pt]{article}
\usepackage{amsmath}
\usepackage{wasysym}
\DeclareMathOperator{\Atang}{Atang.}
\begin{document}
\title{On a new type of rational and highly convergent series, by which
the ratio of the circumference to the diameter is able to be expressed\footnote{Delivered
to the St.--Petersburg Academy June 17, 1779. Originally published as
{\em De novo genere serierum rationalium et valde convergentium, quibus ratio
peripheriae ad diametrum exprimi potest},
Nova Acta Academiae Scientarum Imperialis Petropolitinae \textbf{11} (1798),
150--154, and
republished in \emph{Leonhard Euler, Opera Omnia}, Series 1:
Opera mathematica,
Volume 16, Birkh\"auser, 1992. A copy of the original text is available
electronically at the Euler Archive, at http://www.eulerarchive.org. This paper
is E706 in the Enestr\"om index.}}
\author{Leonhard Euler\footnote{Date of translation: August 8, 2005.
Translated from the Latin
by Jordan Bell, 3rd year undergraduate in Honours Mathematics, School of Mathematics and Statistics, Carleton University,
Ottawa, Ontario, Canada.
Email: jbell3@connect.carleton.ca.
This translation was written
during an NSERC USRA supervised by Dr. B. Stevens.
}}
\date{}
\maketitle

1. The principle, from which these series are deduced, rests in this
binomial formula: $4+x^4$, which is evidently composed of these two
rational factors: $2+2x+xx$ and $2-2x+xx$. Then indeed it at once follows
for this integral formula: $\int \frac{\partial x(2+2x+xx)}{4+x^4}$,
which we shall indicate with the sign $\astrosun$,
to be reduced to this: $\astrosun=\int \frac{\partial x}{2-2x+xx}$, whose
integral, having been obtained so that it vanishes when it is put
$x=0$, is $\Atang \frac{x}{2-x}$.
Whereby it may be observed in the case $x=1$ to be $\astrosun=\frac{\pi}{4}$;
while indeed in the case $x=\frac{1}{2}$ it will be $\astrosun=\Atang \frac{1}{3}$;
then indeed in the case $x=\frac{1}{4}$ it will be $\astrosun=\Atang \frac{1}{7}$.
It is noted moreover for it to be
\[
2\Atang \frac{1}{3}+\Atang \frac{1}{7}=\Atang 1=\frac{\pi}{4}.
\]

2. With therefore this integral formula indicated by the sign $\astrosun$
which is comprised by three parts, each of which we shall unfold separately,
which by the grace of brevity we shall indicate by the following characters:
\[ 
\textrm{I.} \int \frac{\partial x}{4+x^4}=\saturn; \quad
\textrm{II.} \int \frac{x \partial x}{4+x^4}=\jupiter; \quad
\textrm{III.} \int \frac{xx \partial x}{4+x^4}=\mars;
\]
so that it will thus be
\[
\astrosun=2\saturn+2\jupiter+\mars=\Atang \frac{x}{2-x}.
\]
Now therefore we may unfold these three integral formulas in the usual
manner into infinite series, which are thereupon to be formed,
insofar as it will be
\[
\frac{1}{4+x^4}=\frac{1}{4}\Big(1-\frac{x^4}{4}+\frac{x^8}{4^2}-\frac{x^{12}}{4^3}
+
\frac{x^{16}}{4^4}-\textrm{ etc.}\Big).
\]

3. But if now first we adjoin this series with $\partial x$ and we then integrate,
the first formula $\saturn$ will be expressed by the following series:
\[
\saturn=\frac{x}{4}[1-\frac{1}{5}\cdot \frac{x^4}{4}+\frac{1}{9}\Big(\frac{x^4}{4}\Big)^2-
\frac{1}{3}\Big(\frac{x^4}{4}\Big)^3+\textrm{ etc.}]. 
\]
While indeed adjoining the former series by $x\partial x$
and integrating gives
\[
\jupiter=\frac{xx}{8}[1-\frac{1}{3}\cdot \frac{x^4}{4}+\frac{1}{5}\Big(\frac{x^4}{4}\Big)^2-
\frac{1}{7}\Big(\frac{x^4}{4}\Big)^3+\textrm{ etc.}].
\]
Then adjoining the very same series with $xx\partial x$ and integrating
produces
\[
\mars=\frac{x^3}{4}[\frac{1}{3}-\frac{1}{7}\frac{x^4}{4}+
\frac{1}{11}\Big(\frac{x^4}{4}\Big)^2-
\frac{1}{15}\Big(\frac{x^4}{4}\Big)^3+\textrm{ etc.}].
\]

4. With therefore it being $\astrosun=2\saturn+2\jupiter+\mars$, we
shall unfold some particular cases recalled from before, in which it is
$x=1$, $x=\frac{1}{2}$ and $x=\frac{1}{4}$, of which the first is
$\frac{x^4}{4}=\frac{1}{4}$;
for the second indeed it is $\frac{x^4}{4}=\frac{1}{64}$; for the third indeed
$\frac{x^4}{4}=\frac{1}{1024}$; from which it stands open for the two last
cases to converge most greatly, but that the first, whose terms
 decrease by a ratio of four, indeed converges more so than the series
of Leibnitz, by taking an arc whose tangent is $\frac{1}{\sqrt{3}}$,
seeing that this calculation is perturbed by no irrational.

\begin{center}
{\Large The expansion of the first case,\\
where $x=1$ and $\astrosun=\Atang \frac{\pi}{4}$.}
\end{center}

5. Seeing therefore here that it is $\frac{x^4}{4}=\frac{1}{4}$, our three
principal series for $\saturn, \jupiter, \astrosun$ [sic] proceed
in the following way:

\begin{tabular}{l}
$\saturn=\frac{1}{4}[1-\frac{1}{5}\cdot \frac{1}{4}+\frac{1}{9}(\frac{1}{4})^2-
\frac{1}{13}(\frac{1}{4})^3+\frac{1}{17}(\frac{1}{4})^4-\textrm{ etc.}]$\\
$\jupiter=\frac{1}{8}[1-\frac{1}{3}\cdot \frac{1}{4}+\frac{1}{5}(\frac{1}{4})^2-
\frac{1}{7}(\frac{1}{4})^3+\frac{1}{9}(\frac{1}{4})^4-\textrm{ etc.}]$\\
$\mars=\frac{1}{4}[\frac{1}{3}-\frac{1}{7}\cdot \frac{1}{4}+\frac{1}{11}(\frac{1}{4})^2-
\frac{1}{13}(\frac{1}{4})^3+\frac{1}{19}(\frac{1}{4})^4-\textrm{ etc.}]$
\end{tabular}

6. Seeing therefore that it is $\astrosun=2\saturn+2\jupiter+\mars=\frac{\pi}{4}$,
by multiplying the value of $\pi$ by 4, the following three
series are expressed
\begin{align*}
\pi=\begin{cases}
2(1-\frac{1}{5}\cdot \frac{1}{4}+\frac{1}{9}\cdot \frac{1}{4^2}-
\frac{1}{13}\cdot \frac{1}{4^3}+\frac{1}{17}\cdot \frac{1}{4^4}-
\textrm{ etc.})\\
1(1-\frac{1}{3}\cdot \frac{1}{4}+\frac{1}{5}\cdot \frac{1}{4^2}-
\frac{1}{7}\cdot \frac{1}{4^3}+\frac{1}{9}\cdot \frac{1}{4^4}-\textrm{ etc.})\\
1(\frac{1}{3}-\frac{1}{7}\cdot \frac{1}{4}+\frac{1}{11}\cdot \frac{1}{4^2}-
\frac{1}{15}\cdot \frac{1}{4^3}+\frac{1}{19}\cdot \frac{1}{4^4}-
\textrm{ etc.})
\end{cases}
\end{align*}

7. From these particular three series the ratio of the circumference to
the diameter is able to be calculated with much less work
than by the series of Leibnitz, which method the most meritous Authors
Sharp, Machin and de Lagny have used, of whom the first has determined
$\pi$ in a decimal fraction to 72 figures, the second to 100, and the last
indeed to 128. And truly they were able to lift up the following cases with
much more effort.

\begin{center}
{\Large The expansion of the second case,\\
where $x=\frac{1}{2}$.}
\end{center}

8. In this case it will therefore be $\frac{x^4}{4}=\frac{1}{64}$, from
which the three series are drawn forth in the following way:

\begin{tabular}{l}
$\saturn=\frac{1}{8}(1-\frac{1}{5}\cdot \frac{1}{64}+\frac{1}{9}\cdot
\frac{1}{64^2}-\frac{1}{13}\cdot \frac{1}{64^3}+\textrm{ etc.})$\\
$\jupiter=\frac{1}{32}(1-\frac{1}{3}\cdot \frac{1}{64}+\frac{1}{5}\cdot
\frac{1}{64^2}-\frac{1}{7}\cdot \frac{1}{64^3}+\textrm{ etc.})$\\
$\mars=\frac{1}{32}(\frac{1}{3}-\frac{1}{7}\cdot \frac{1}{64}+\frac{1}{11}
\cdot \frac{1}{64^2}-\frac{1}{15}\cdot \frac{1}{64^3}+\textrm{ etc.})$
\end{tabular}

9. Therefore with it $2\saturn+2\jupiter+\mars=\Atang \frac{1}{3}$, it will
be
\begin{align*}
\Atang \frac{1}{3}=\begin{cases}
\frac{1}{4}(1-\frac{1}{5}\cdot \frac{1}{64}+\frac{1}{9}\cdot \frac{1}{64^2}-
\frac{1}{13}\cdot \frac{1}{64^3}+\textrm{ etc.})\\
\frac{1}{16}(1-\frac{1}{3}\cdot \frac{1}{64}+\frac{1}{5}\cdot \frac{1}{64^2}-
\frac{1}{7}\cdot \frac{1}{64^3}+\textrm{ etc.})\\
\frac{1}{32}(\frac{1}{3}-\frac{1}{7}\cdot \frac{1}{64}+\frac{1}{11}\cdot
\frac{1}{64^2}-\frac{1}{15}\cdot \frac{1}{64^3}+\textrm{ etc.})
\end{cases}
\end{align*}
Even though here these three series are to be computed, however,
because each successively decreases by the same ratio $1:64$, this
labor will be able to be shortened in a wonderful way.

\begin{center}
{\Large The expansion of the third case,\\
where $x=\frac{1}{4}$.}
\end{center}

10. Seeing therefore here that it is $\frac{x^4}{4}=\frac{1}{1024}$,
our three principal series will be had as follows:

\begin{tabular}{l}
$\saturn=\frac{1}{16}(1-\frac{1}{5}\cdot \frac{1}{1024}+\frac{1}{9}\cdot
\frac{1}{1024^2}-\frac{1}{13}\cdot \frac{1}{1024^3}+\textrm{ etc.})$\\
$\jupiter=\frac{1}{128}(1-\frac{1}{3}\cdot \frac{1}{1024}+
\frac{1}{5}\cdot \frac{1}{1024^2}-\frac{1}{7}\cdot \frac{1}{1024^3}+\textrm{ etc.})$\\
$\mars=\frac{1}{256}(\frac{1}{3}-\frac{1}{7}\cdot \frac{1}{1024}+
\frac{1}{11}\cdot \frac{1}{1024^2}-\frac{1}{15}\cdot \frac{1}{1024^3}+
\textrm{ etc.})$
\end{tabular}

11. Therefore with $2\saturn+2\jupiter+\mars=\Atang \frac{1}{7}$, it will
properly be by joining these series:
\begin{align*}
\Atang \frac{1}{7}=\begin{cases}
\frac{1}{8}(1-\frac{1}{5}\cdot \frac{1}{1024}+\frac{1}{9}\cdot \frac{1}{1024^2}-
\frac{1}{13}\cdot \frac{1}{1024^3}+\textrm{ etc.})\\
\frac{1}{64}(1-\frac{1}{3}\cdot \frac{1}{1024}+\frac{1}{5}\cdot \frac{1}{1024^2}-
\frac{1}{7}\cdot \frac{1}{1024^3}+\textrm{ etc.})\\
\frac{1}{256}(\frac{1}{3}-\frac{1}{7}\cdot \frac{1}{1024}+
\frac{1}{11}\cdot \frac{1}{1024^2}-
\frac{1}{15}\cdot \frac{1}{1024^3}+\textrm{ etc.})
\end{cases}
\end{align*}

\begin{center}
{\Large An application of the last two cases for expressing a circumference
 by highly convergent series.}
\end{center}

12. With, as we have already observed, 
$\frac{\pi}{4}=2\Atang \frac{1}{3}+\Atang \frac{1}{7}$ it will be
$\pi=8\Atang \frac{1}{3}+4\Atang \frac{1}{7}$,
by substituting the value of $\pi$ into the series found above,
the following six series may be expressed together:
\begin{align*}
\pi=\begin{cases}
2(1-\frac{1}{5}\cdot \frac{1}{64}+\frac{1}{9}\cdot \frac{1}{64^2}-
\frac{1}{13}\cdot \frac{1}{64^3}+\textrm{ etc.})\\
\frac{1}{2}(1-\frac{1}{3}\cdot \frac{1}{64}+\frac{1}{5}\cdot \frac{1}{64^2}-
\frac{1}{7}\cdot \frac{1}{64^3}+\textrm{ etc.})\\
\frac{1}{4}(\frac{1}{3}-\frac{1}{7}\cdot \frac{1}{64}+\frac{1}{11}\cdot \frac{1}{64^2}-
\frac{1}{15}\cdot \frac{1}{64^3}+\textrm{ etc.})\\
\frac{1}{2}(1-\frac{1}{5}\cdot \frac{1}{1024}+\frac{1}{9}\cdot \frac{1}{1024^2}-
\frac{1}{13}\cdot \frac{1}{1024^3}+\textrm{ etc.})\\
\frac{1}{16}(1-\frac{1}{3}\cdot \frac{1}{1024}+\frac{1}{5}\cdot \frac{1}{1024^2}-
\frac{1}{7}\cdot \frac{1}{1024^3}+\textrm{ etc.})\\
\frac{1}{64}(1-\frac{1}{7}\cdot \frac{1}{1024}+\frac{1}{11}\cdot \frac{1}{1024^2}-
\frac{1}{15}\cdot \frac{1}{1024^3}+\textrm{ etc.})
\end{cases}
\end{align*} 
This occurs as most noteworthy, because all these series proceed only by
powers of two.

\end{document}